\documentclass{compositio}

\usepackage{amssymb, mathrsfs, amsmath}
\usepackage[all]{xy}

\newcommand{\C}{\mathcal{C}}

\newcommand{\E}{\mathcal{E}}

\newcommand{\HH}{\mathbb{H}}

\newcommand{\OO}{\mathcal{O}}

\newcommand{\ZZ}{\mathbb{Z}}
\newcommand{\zll}{{\mathbb{Z}_{(l)}}}
\newcommand{\zpi}{\mathbb{Z}[\tfrac{1}{p}]}

\newcommand{\id}{\mathrm{id}}
\newcommand{\Tr}{\mathrm{Tr}}
\newcommand{\Perf}{\mathrm{Perf}}

\newtheorem{theo}{Theorem}[section]
\newtheorem{theoUn}[theo]{Theorem}
\newtheorem{coro}[theo]{Corollary}
\newtheorem{lemm}[theo]{Lemma}
\newtheorem{prop}[theo]{Proposition}
\theoremstyle{definition}
\newtheorem{defi}[theo]{Definition}
\theoremstyle{remark}
\newtheorem{rema}[theo]{Remark}

\newcommand{\fpsl}{fps$l'$}
\newcommand{\fpslgl}{fps$l'$gl}
\newcommand{\ldh}{{$l$dh}}

\title{Vanishing of negative $K$-theory in positive characteristic}
\author{Shane Kelly}
\address{Interactive Research Center of Science\\
Graduate School of Science and Engineering\\
Tokyo Institute of Technology\\
2-12-1 Ookayama, Meguro\\
Tokyo 152-8551 JAPAN }
\thanks{This work was done as part of the authors Ph.D., during which he was supported by an Australian National University Vice-Chancellor's Scholarship and an Australian National University PhD Supplementary Scholarship.}

\begin{document}

\begin{abstract}
We show how a theorem of Gabber on alterations can be used to apply work of Cisinski, Suslin, Voevodsky, and Weibel to prove that $K_n(X) \otimes \zpi = 0$ for $n < - \dim X$ where $X$ is a quasi-excellent noetherian scheme, $p$ is a prime that is nilpotent on $X$, and $K_n$ is the $K$-theory of Bass-Thomason-Trobaugh. This gives a partial answer to a question of Weibel.
\end{abstract}

\maketitle

\section{Introduction}

In \cite[2.9]{Wei80} Weibel asks if $K_n(X) = 0$ for $n < - \dim X$ for every noetherian scheme $X$ where $K_n$ is the $K$-theory of Bass-Thomason-Trobaugh. This question was answered in the affirmative in \cite{CHSW} for schemes essentially of finite type over a field of characteristic zero. Assuming strong resolution of singularities, it is also answered in the affirmative in \cite{GH10} for schemes essentially of finite type over an infinite perfect field of positive characteristic. Both of these proofs compare $K$-theory with cyclic homology, and then use a cdh descent argument.

Our main theorem is the following.

\begin{theoUn}[{\ref{theo:KtheoryVanishing}}]
Let $X$ be a quasi-excellent noetherian scheme and $p$ a prime that is nilpotent on $X$. Then $K_n(X) \otimes \zpi = 0$ for $n < - \dim X$ where $K_n$ is the $K$-theory of Bass-Thomason-Trobaugh.
\end{theoUn}

Our proof can be outlined as follows. We reduce the vanishing of $K_n(X) \otimes \zpi$ to the vanishing of homotopy invariant $K$-theory using a result of Weibel that in positive characteristic $p$ compares $K$-theory with homotopy invariant $K$-theory $KH_n$ away from $p$-torsion. Cisinski has shown that homotopy invariant $K$-theory is representable in the Morel-Voevodsky stable homotopy category \cite[Theorem 2.20]{Cis13}. Furthermore, he shows using Ayoub's proper base change theorem (\cite[Corollary 1.7.18]{Ayo07}) and Voevodsky's work on cd structures (\cite{Voe10a}, \cite{Voe10b}) that every object in the stable homotopy category satisfies cdh descent \cite[Proposition 3.7]{Cis13}. In particular, homotopy invariant $K$-theory satisfies cdh descent for noetherian schemes of finite Krull dimension (without restriction on the characteristic).\footnote{That homotopy invariant $K$-theory satisfies cdh descent was shown in characteristic zero in \cite{Hae04} for schemes essentially of finite type over a field and the proof goes through in positive characteristic if embedded resolution of singularities turns out to hold. See also \cite[Theorem 5.5]{Wei81} for a cdh descent result for $KH \otimes [\tfrac{1}{p}]$ over $p$-torsion bases.} Using the cdh descent spectral sequence, the desired vanishing is implied by the vanishing of the cdh sheaves associated to $KH_n \otimes \zpi$ for $n < 0$. To deduce this, we apply a theorem of Gabber on alterations as a replacement for resolution of singularities.

The cdh topology was introduced to apply resolution of singularities. To apply this theorem of Gabber, we refine the cdh topology to a topology that we call the \ldh~topology. This topology is by definition generated by the cdh topology, and morphisms which are flat finite surjective and of constant degree prime to a chosen prime $l$. This method of applying the \ldh~topology to replace resolution of singularities arguments with Gabber's theorem on alterations can also be used in other settings. See for example \cite{KellyThesis} where all the results of \cite{Voev00}, \cite{FV}, and \cite{Sus00} that assume resolution of singularities are proved in positive characteristic using this theorem of Gabber (using $\zpi$-coefficients of course, where $p$ is the characteristic of the base field). See also \cite{HKO} where the \ldh~topology is used to generalise to positive characteristic Voevodsky's result identifying the algebra of bistable operations with the motivic Steenrod algebra.

\acknowledgements{The material in this article is part of my Ph.D. thesis. I am deeply indebted to Denis-Charles Cisinski for many discussions, for the choice of problem, and in particular for the myriad potential (not all equivalent) definitions for the \ldh~topology that he suggested. My original proof was much longer, and this much shorter one was discovered while preparing for a talk I gave in the Oberseminar in Algebraic Geometry at The University of Zurich. I thank them for the invitation. I also thank Thomas Geisser and Chuck Weibel for pointing out missing references. The readability of this paper was greatly improved following suggestions by Ariyan Javanpeykar, Giuseppe Ancona, and the referees.}

\section{The \ldh~topology} \label{sec:ltopology}

We work with $Sch(S)$, the category of separated schemes of finite type over a noetherian base scheme $S$. Recall that a refinement of a family of morphisms $\{ U_i \to X \}_{i \in I}$ is a family of morphisms $\{ V_j \to X \}_{j \in J}$ such that for each $j \in J$ there is an $i_j \in I$ and a factorisation $V_j \to U_{i_j} \to X$. Following Suslin and Voevodsky (and contrary to Artin and therefore Milne), we use the terms \emph{topology} and \emph{covering} as in \cite[Definition 1.1]{SGA4II} and \cite[Definition 1.2]{SGA4II} respectively. In particular, if a family of morphisms admits a refinement by a covering for a topology $\tau$, then that family itself is also a covering \cite[Proposition 1.4]{SGA4II}. We observe the usual abuse of the term \emph{covering} by referring to a morphism $f$ as a covering if $\{ f \}$ is a covering family.

The reader not familiar with the cdh topology can find a very readable account of it in \cite[Section 5]{SV00}.

\begin{samepage}
\begin{defi} \label{defi:topologies}
Let $l \in \ZZ$ be a prime.
\begin{enumerate}
 \item By an \emph{\fpsl~morphism} (fini, plat, surjectif, premier {\`a} $l$) we will mean a morphism $f: U \to X$ that is finite flat surjective of constant rank prime to $l$.
 \item The \emph{\fpsl~topology} on the category $Sch(S)$ is the least fine topology such that \fpsl~morphisms are \fpsl~coverings.
 \item The \emph{\ldh~topology} on the category $Sch(S)$ is the least fine topology which is finer than the cdh topology, and finer than the \fpsl~topology.
\end{enumerate}
\end{defi}
\end{samepage}

Our choice of definition of an \ldh~topology is one of the reasons our vanishing argument is so direct. It is motivated by the following two ideas. Firstly, the theorem of Gabber (Theorem~\ref{theo:gabber}) should provide the existence of regular \ldh~coverings (or smooth depending on the context). Secondly, we want to make use of the ample literature on the cdh topology. That is, we want to be able to reduce statements about the \ldh~topology to statements about the cdh topology and statements about the \fpsl~topology. This way we only need to deal with the \fpsl~topology. This we usually do using a suitable structure of ``trace'' morphisms.

\begin{rema}
Our \ldh~topology differs from the topology of $\ell'$ alterations\footnote{We have chosen to use $l$ instead of the more aesthetically pleasing $\ell$ used in \cite{ILO} as $l$ is easier for search engines.} described in \cite[Expos{\'e} III, Section 3]{ILO}.

Firstly, the underlying categories are different. In \cite{ILO} the category $alt/S$ is used, which should be thought of as some kind of Riemann-Zariski space. The objects of this category are the reduced schemes that are dominant and of finite type over $S$, such that the structural morphism sends generic points to generic points, and the induced field extensions associated to generic points are finite. The inclusion $alt/S \to Sch(S)$ does not preserve fibre products.\footnote{The category $alt/S$ has fibre products in the categorical sense. These are obtained from the usual fibre product in the category of schemes as follows. For morphisms $X \to Y, Z \to Y$ in $alt/S$ one forms the usual fibre product $X \times_Z Y$ in the category of schemes, then takes the associated reduced scheme $(X \times_Z Y)_{red}$, and omits any irreducible components of $(X \times_Z Y)_{red}$ which do not dominate an irreducible component of $S$ \cite[Expos{\'e} II, 1.2.5, Proposition 1.2.6]{ILO}.} If we equip $alt/S$ with the topology of $\ell'$ alterations and $Sch(S)$ with the \ldh~topology, then this inclusion is not a continuous morphism of sites \cite[D{\'e}finition 1.1]{SGA4III} (it is however cocontinuous \cite[D{\'e}finition 2.1]{SGA4III}).

Second, the topology of $\ell'$ alterations is defined in \cite{ILO} similarly to the ($alt/S$  analogue of the) cdh topology, but by extending the ``proper cdh'' part to allow Gabber's theorem to replace resolution of singularities. It appears that our topology is a ``global'' version of their ``local'' topology where global and local are in the resolution of singularities sense (cf. \cite[Expos{\'e} II, Th{\'e}or{\`e}me 3.2.3]{ILO}). 
\end{rema}

The next lemma and proposition are of the form: ``a $\sigma$ covering followed by a $\rho$ covering can be refined by a $\rho$ covering followed by a $\sigma$ covering''. After that, Lemma~\ref{lemm:sep} says that in this situation, the associated $\sigma$ sheaf functor preserves $\rho$ separated presheaves.

\begin{lemm} \label{lemm:NisFpslSwap}
Suppose that $\{Y_i' \to Y\}_{i \in I}$ is a Nisnevich covering and $Y \to X$ a finite morphism. Then there exists a Nisnevich covering $\{X_j' \to X\}_{j \in J}$ such that $\{W_{jk} \to Y\}_{j \in J, k \in K_j}$ is a refinement of $\{Y_i' \to Y\}_{i \in I}$ where the $W_{jk}$ are the connected components of $Y \times_X X_j'$.


\end{lemm}

\begin{proof}
If $X$ is henselian, then $Y$ is also henselian \cite[Proposition 18.5.10]{EGAIV4} and so there exists some $i \in I$ for which $Y_i' \to Y$ has a section. So we can take $\{X_j' \to X\}_{j \in J}$ to be $\{X \stackrel{id}{\to} X \}$. If $X$ is not henselian then for every point $x \in X$ we consider the pullback along the henselisation $\ ^hx \to X$. The result now follows from the limit arguments in \cite[Section 8]{EGAIV3} and the description of the henselisation as a suitable limit of \'{e}tale neighbourhoods \cite[Definition 18.6.5]{EGAIV3}.
\end{proof}

The following proposition is in a similar spirit to \cite[Proposition 5.9]{SV00}.

\begin{prop} \label{prop:cdhlswap}
Let $X$ be a noetherian scheme and suppose that $Y \to X$ is an \fpsl~morphism and $\{ U_i \to Y \}_{i \in I}$ is a cdh covering. Then there exists a cdh covering $\{ V_j \to X \}_{j \in J}$ and \fpsl~morphisms $V'_j \to V_j$ such that $\{ V'_j \to X \}_{j \in J}$ is a refinement of $\{ U_i \to Y \to X \}_{i \in I}$.
\end{prop}

\begin{proof}
The proof is similar to the proof of \cite[Proposition 5.9]{SV00}, and the reader is encouraged to consult that proof before proceeding.

We first consider the case when $\{ U_i \to Y \}_{i \in I}$ is a singleton $\{U \to Y\}$ containing a proper morphism which is a cdh covering (such a morphism is called a \emph{proper cdh covering}). As in the proof of \cite[Proposition 5.9]{SV00}, proceeding by noetherian induction we may assume that for any proper closed subscheme $Z \subset X$ the induced covering of $Z$ admits a refinement of the desired form and the scheme $X$ is integral.

Suppose we can find a proper birational morphism $V \to X$ and an \fpsl~covering $V' \to V$ such that the composition $V' \to V \to X$ factors through the composition $U \to Y \to X$. Choose a closed subscheme $Z \to X$ outside of which $V \to X$ is an isomorphism. The inductive hypothesis applied to $\{Z \times_X U \to Z \times_X Y \}$ and $Z \times_X Y \to Z$ over $Z$ gives a cdh covering $\{W_j \to Z \}_{j \in J}$ and \fpsl~morphisms $W'_j \to W_j$ such that $\{ W'_j \to W_j \to Z \}_{j \in J}$ is a refinement of $\{Z \times_X U \to Z \times_X Y \to Z\}$. Consequently, $\{W'_j \to W_j \to Z \to X \}_{j \in J} \cup \{V' \to V \to X\}$ is a refinement of $\{U \to Y \to X \}$ of the desired form. So it suffices to find such a $V \to X$ and $V' \to V$.

We will now construct the following diagram, which is presented here so that the reader may refer to it while it is being constructed (we assure the anxious reader that the symbols are defined in the following paragraph, during the construction of the diagram).
\[ \xymatrix{
V' \ar[r] \ar[dd] & Y'' \times_Y U \ar[rr] \ar[d] && U \ar[d] \\
& Y'' \ar[r] \ar[d] & \overline{\{ \eta_i \}} \ar[d]  \ar[r] & Y \ar[dl] \\
V \ar[r] & X' \ar[r] & X
} \]
Replacing $U \to Y$ by a finer proper cdh covering we may assume that $U_{red} \to Y_{red}$ is an isomorphism over a dense open subscheme of $Y$, and that $U = U_{red}$. If $\eta_i$ are the generic points of $Y$ and $m_i$ the lengths of their local rings, then the degree of $Y \to X$ is $\sum m_i [k(\eta_i): k(\xi)]$ where $\xi$ is the generic point of $X$. Since $l$ does not divide $\sum m_i [k(\eta_i): k(\xi)]$, there is some $i$ for which $l$ does not divide $[k(\eta_i): k(\xi)]$. Choose such an $i$ and consider the closed integral subscheme $\overline{\{ \eta_i \}}$ of $Y$ which has $\eta_i$ as its generic point. By the platification theorem \cite{RG71} (or \cite[Theorem 2.2.2]{SV} for a precise statement) there exists a blow-up $X' \to X$ of $X$ with nowhere dense centre such that the strict transform $Y'' \to X'$ of $\overline{\{ \eta_i \}} \to X$ is flat, and hence finite flat surjective of degree prime to $l$. Consequently, the composition $Y'' \times_Y U \to Y'' \to X'$ is generically an \fpsl~covering. Applying the platification theorem to this composition we can find a blow-up $V \to X'$ of $X'$ such that the strict transform $V' \to V$ of the composition $Y'' \times_Y U \to Y'' \to X'$ is flat, Since it is generically an \fpsl~covering, flatness implies that it is actually an \fpsl~covering.

Having concluded the case when $\{ U_i \to Y \}_{i \in I}$ is a singleton containing a proper cdh covering, we now return the general case. Using \cite[Proposition 5.9]{SV00}, we can assume that $\{ U_i \to Y \}_{i \in I}$ is of the form $\{U_i \to Y' \to Y\}^n_{i = 1}$ where $Y' \to Y$ is a proper cdh covering and $\{U_i \to Y' \}_{i = 1}^n$ is a Nisnevich covering of $Y'$. The proper cdh case that we just considered provides a cdh covering $\{ V_j \to X \}_{j \in J}$ and \fpsl~morphisms $V'_j \to V_j$ such that $\{ V'_j \to X \}_{j \in J}$ is a refinement of $\{ Y' \to Y \to X \}$. Now for each $j \in J$ apply Lemma~\ref{lemm:NisFpslSwap} to the Nisnevich covering $\{V'_j \times_{Y'} U_i \to V'_j\}_{i \in I}$ and the finite morphism $V'_j \to V_j$.
\end{proof}

For a topology $\tau$ we denote by $a_\tau$ the canonical functor that takes a presheaf to its associated $\tau$ sheaf.

\begin{lemm} \label{lemm:sep}
Suppose that $\sigma$ and $\rho$ are two topologies on a category $\C$. Suppose that for every $\rho$ covering $\{U_i \to X \}_{i \in I}$ and $\sigma$ coverings $\{ V_{ij} \to U_i\}_{j \in J_i}$, there exists a $\sigma$ covering $\{U_k' \to X \}_{k \in K}$ and $\rho$ coverings $\{ V_{k\ell}' \to U_k'\}_{\ell \in L_k}$ such that $\{ V_{k\ell}' \to X\}_{k \in K, \ell \in L_k}$ refines $\{ V_{ij} \to X\}_{i \in I, j \in J_i}$. If $F$ is a presheaf which is $\rho$ separated, then its associated $\sigma$ sheaf $a_\sigma F$ is also $\rho$ separated.
\end{lemm}

\begin{rema}
The family $\mathcal{V} = \{ V_{k\ell}' \to X\}_{k \in K, \ell \in L_k}$ is not required to be a covering for either of the topologies $\sigma$ or $\rho$. In practice, this lemma will be applied in situations where we are interested in the topology generated by $\sigma$ and $\rho$ (the family $\mathcal{V}$ is a covering for this topology).
\end{rema}

\begin{proof}
We will use coverings of cardinality one to make the proof easier to read. Suppose $s \in a_\sigma F(X)$ is a section and $U \to X$ a $\rho$ covering such that $s|_U = 0$. We must show $s = 0$ in $a_\sigma F(X)$. It is sufficient to consider the case where $s$ is in the image of $F(X) \to a_\sigma F(X)$. In this case the condition $s|_U = 0$ implies that there is a $\sigma$ covering $V \to U$ with $s'|_V = 0$ where $s' \in F(X)$ is some section sent to $s \in a_\sigma F(X)$. By hypothesis, we can refine $V \to U \to X$ by a composition $V' \to U' \to X$ with $V' \to U'$ a $\rho$ covering and $U' \to X$ a $\sigma$ covering. Now $s'|_{V'} = 0$ but $F$ is $\rho$ separated so $s'|_{U'} = 0$, which implies that $s = 0$ in $a_\sigma F(X)$.
\end{proof}

We now reproduce a weak version of a theorem of Gabber. We follow it with a corollary which converts it into a form that we will use. For a statement and an outline of the proof of this theorem of Gabber see \cite{Ill09}, or \cite{Gab05}. There is a proof in the book in preparation \cite{ILO}.


\begin{theo}[{(Gabber), \cite[Expos{\'e} 0 Theorem 2]{ILO}, \cite[Expos{\'e} II Theorem 3.2.1]{ILO}}] \label{theo:gabber}
Let $X$ be a noetherian quasi-excellent scheme and let $l$ be a prime number invertible on $X$. There exists a finite family of morphisms of finite type $\{U_i \to X \}_{i = 1}^n$ with each $U_i$ regular, and a refinement of this family of the form $\{ V_j \to Y \to X \}_{j \in J}$ such that $Y \to X$ is proper surjective generically finite of degree prime to $l$ and $\{V_j \to Y \}_{j \in J}$ is a Nisnevich covering.
\end{theo}

\begin{coro} \label{coro:regularlCover}
Let $X$ be a noetherian quasi-excellent scheme and let $l$ be a prime number invertible on $X$. Then there exists an \ldh~covering $\{ W_i \to X \}_{i = 1}^m$ of $X$ such that each $W_i$ is regular.
\end{coro}

\begin{proof}
This corollary follows from Theorem~\ref{theo:gabber} using the platification theorem \cite{RG71}.

Let us be more explicit. We will use noetherian induction. Suppose that the result is true for all proper closed subschemes of $X$. Let $\{U_i \to X \}_{i = 1}^n$ and $\{ V_j \to Y \to X \}_{j \in J}$ be as in the statement of Theorem~\ref{theo:gabber}. As in the proof of Proposition~\ref{prop:cdhlswap} (or rather, the proof of \cite[Proposition 5.9]{SV00}) we can assume that $X$ is integral, since the disjoint union of the inclusion of the reduced irreducible components of $X$ is a cdh covering, and hence an \ldh~covering. By the platification theorem there exists a blow-up with nowhere dense centre $X' \to X$ such that the proper transform $Y' \to X'$ of $Y \to X$ is a finite flat surjective morphism of constant degree. Let $\{ V'_j \to Y' \}_{j \in J}$ be the pullback of the Nisnevich covering $\{ V_j \to Y \}_{j \in J}$, and let $Z \subset X$ be a proper closed subscheme of $X$ for which $X' \to X$ is an isomorphism outside of $Z$. The family $\{Z \to X\} \cup \{V'_j \to X \}_{j \in J}$ is the composition of the cdh covering $\{Z \to X, X' \to X\}$, the \fpsl~covering $\{Z \stackrel{id}{\to} Z, Y' \to X' \}$, and the Nisnevich covering $\{Z \stackrel{id}{\to} Z\} \cup \{V'_j \to Y' \}_{j \in J}$. Hence, it is an \ldh~covering. Now by construction this family is a refinement of the family $\{Z \to X\} \cup \{U_i \to X \}_{i = 1}^n$. Therefore $\{Z \to X\} \cup \{U_i \to X \}_{i = 1}^n$ is also an \ldh~covering. Finally, by the inductive hypothesis there exists an \ldh~covering $\{Z_k' \to Z\}_{k= 1}^{n'}$ of $Z$ with each $Z_k'$ regular, and so the family $\{Z_k' \to Z \to X\}_{k= 1}^{n'} \cup \{U_i \to X \}_{i = 1}^n$ is an \ldh~covering of $X$ for which all the sources are regular.
\end{proof}

\section{Vanishing of $K$-theory} \label{sec:vanishing}

In this section we will use the $K$-theory of Bass-Thomason-Trobaugh and Weibel's homotopy invariant $K$-theory. These will be denoted by $K$ and $KH$ respectively. 

Let us quickly recall their constructions. For any scheme $X$ consider $\Perf(X)$, the complicial biWaldhausen category of perfect complexes of $\OO_X$-modules \cite[Definition 1.2.11, Definition 2.2.12]{TT90}. We deviate slightly from \cite{Wei89} and \cite{TT90} by using $\OO_X$-modules on $Sch(X)$ as opposed to the small Zariski site of $X$. We do this so that $K$ and $KH$ are actual presheaves instead of just lax functors. See \cite[Section C4]{FS02} for a discussion. From the biWaldhausen category $\Perf(X)$ we obtain a spectrum denoted in \cite[6.4]{TT90} by $K^B(X)$, but which will be denoted in our article simply by $K(X)$. The $K(X)$ form a presheaf of spectra $K$ on $Sch(X)$. From this presheaf of spectra we obtain a second presheaf of spectra \cite[Definition 1.1]{Wei89} which is denoted by $KH$. There is a canonical morphism of presheaves of spectra $K \to KH$.

We write $KH_n$ (resp. $K_n$) for the functor which associates to a scheme $X$ the $n$th homotopy group of $KH(X)$ (resp. $K(X)$).

It is in the following lemma that we use ``trace'' morphisms to bridge the gap between cdh and \ldh.

\begin{lemm} \label{lemm:KHfpslsep}
The cdh sheaves $a_{cdh}K_n \otimes \zll$ associated to the presheaves $K_n \otimes \zll$ are separated for the \ldh~topology.
\end{lemm}

\begin{proof}
Consider the weakest topology on $Sch(S)$ such that \fpsl~morphisms $f: Y \to X$ with $f_*\OO_Y$ globally free are coverings. Let's call this topology the \fpslgl~topology. We begin by showing that the presheaves $K_n \otimes \zll$ are separated for the \fpslgl~topology.

The construction of $K$ is functorial in complicial biWaldhausen categories. Notably, for every finite flat surjective morphism of schemes $f: Y \to X$ we get a corresponding exact functor $f_*: \Perf(Y) \to \Perf(X)$ between the biWaldhausen categories. These functors give rise to morphisms $\Tr_{f} : K_n(Y) \to K_n(X)$. We claim that if there is an isomorphism $f_*\OO_Y \cong \OO_X^d$ then $\Tr_{f} K_n(f) = d \cdot id_{K_n(X)}$. If $f$ is an \fpsl~morphism then $d$ is invertible in $\zll$ and we deduce from $\Tr_{f} K_n(f) = d \cdot id_{K_n(X)}$ that the morphism $K_n(f) \otimes \zll$ is injective, and hence $K_n \otimes \zll$ is separated for the \fpslgl~topology. To prove $\Tr_{f} K_n(f) = d \cdot id_{K_n(X)}$, since $f_*\OO_Y$ is globally free, we reduce immediately to proving that the operation of $A \mapsto A^{\oplus d}$ on $\Perf(X)$ induces $a \mapsto d \cdot a$ in $K_n$. This can be derived from the additivity theorem \cite[1.7.3]{TT90} by induction on $d$ (let $F' = \id, F'' = (-)^{\oplus (d - 1)}$, $F = - \oplus (-)^{\oplus (d - 1)}$ with the obvious natural transformations).

Now Lemma~\ref{lemm:NisFpslSwap} shows that the hypotheses of Lemma~\ref{lemm:sep} hold with $\sigma =$ Nis and $\rho =$ \fpslgl. Hence, the Nisnevich sheaf $a_{Nis}K_n \otimes \zll$ associated to $K_n \otimes \zll$ is \fpslgl~separated. But every finite flat morphism of schemes is free locally for the Zariski topology, and therefore locally for the Nisnevich topology. Hence, $a_{Nis}K_n \otimes \zll$ is separated for the \fpsl~topology. Now Proposition~\ref{prop:cdhlswap} shows that the hypotheses of Lemma~\ref{lemm:sep} hold with $\sigma =$ cdh and $\rho =$ \fpsl, so the associated cdh sheaves $a_{cdh}a_{Nis}K_n \otimes \zll$ are \fpsl~separated. Clearly we have $a_{cdh}a_{Nis}K_n \otimes \zll = a_{cdh}K_n \otimes \zll$ as the cdh topology is finer than the Nisnevich topology by definition. Now it follows from the definition of the \ldh~topology that $a_{cdh}K_n \otimes \zll$ is separated for the \ldh~topology.
\end{proof}

\begin{defi}
For $\E$ a presheaf of spectra on $Sch(S)$ and $\tau$ a topology with enough points define $\HH_\tau(-, \E)$ to be the presheaf of spectra given by the Godement-Thomason construction \cite[1.33]{Tho85}. This comes equipped with a canonical morphism $\E \to \HH_\tau(-, \E)$.
\end{defi}

The following lemma shows that we can apply the previous definition to the cdh topology on $Sch(S)$.

\begin{lemm}[(Deligne, Suslin-Voevodsky)]
The cdh topology on $Sch(S)$ has enough points.
\end{lemm}


\begin{proof}
A theorem of Deligne says that any locally coherent topos has enough points \cite[Proposition 9.0]{SGA4VI}.\footnote{The author thanks Brad Drew for pointing out that there is an extremely clear account of this theorem of Deligne given in \cite[Chapter 7]{Joh77}. The statement is \cite[Theorem 7.44]{Joh77}.} A topos is locally coherent \cite[Definition 2.3]{SGA4VI} if and only if it is equivalent to a category of sheaves on a site such that every object is quasi-compact, and all fibre products and finite products are representable \cite[2.4.5]{SGA4VI}. By definition \cite[Definition 1.1]{SGA4VI}, an object is quasi-compact if and only if every covering family admits a finite subfamily which is still a covering family. The proof of \cite[Proposition 5.9]{SV00} remains valid with the base field $F$ replaced by any noetherian scheme $S$, and shows that every object in $Sch(S)$ is quasi-compact (for the cdh topology).
\end{proof}

\begin{theo}[{\cite[Theorem 2.20, Proposition 3.7]{Cis13}}] \label{theo:KHcdhDescent}
The presheaf of spectra $KH$ satisfies cdh descent on $Sch(S)$. That is, the canonical morphism $KH \to \HH_{cdh}(-, KH)$ gives a stable weak equivalence of spectra when evaluated on each scheme.
\end{theo}

\begin{theo} \label{theo:KtheoryVanishing}
Let $X$ be a quasi-excellent noetherian scheme and $p$ a prime that is nilpotent on $X$. Then $K_n(X) \otimes \zpi = 0$ for $n < - \dim X$.
\end{theo}

\begin{proof}
Since $p$ is nilpotent on $X$ the canonical morphism $K_n \otimes \zpi \to KH_n \otimes \zpi$ is an isomorphism \cite[9.6]{TT90}. Hence it suffices to prove that \mbox{$KH_n(X) \otimes \zll$} vanishes for every prime $l \neq p$ and $n < - \dim X$. Since $KH$ satisfies cdh descent (Theorem~\ref{theo:KHcdhDescent}) we have a spectral sequence of the following form\footnote{The indexing we have used is the Bousfield-Kan indexing used in \cite{Tho85}. Other indexing conventions exist in the literature. See \cite[Corollary 5.2]{Wei89} and \cite[Section 5]{GH10} for examples of other possible indexings.} (cf. \cite[1.36]{Tho85}).
\[ E_2^{p, q} = H^{p}_{cdh}(X, a_{cdh}KH_q(-)) \Rightarrow KH_{q - p}(X) \]
As the reader may be unfamiliar with this indexing, we mention that the differentials on the $E_r$ sheet have bidegree $(r, r - 1)$. Due to the cdh cohomological dimension of $X$ being bounded by $\dim X$ \cite[12.5]{SV00}, the $E_2$ sheet is zero outside of $0 \leq p \leq \dim X$. This implies that the spectral sequence converges completely (even strongly) to $KH_n(X)$ for all $n$ \cite[Lemma 5.48]{Tho85}. So the filtration on $KH_n(X)$ is complete and Hausdorff (it is automatically exhaustive), and the filtration quotients are isomorphic to $E_\infty^{p,p + n}$ \cite[Definition 5.46]{Tho85}. As $-\otimes \zll$ is exact, to show the desired vanishing of \mbox{$KH_n(X) \otimes \zll$}, it therefore suffices to show that
\[ E_\infty^{p,p + n} \otimes \zll = 0 \quad  \textrm{ for all } \quad p \in \ZZ, n < \dim X. \]
Since we already know that the $E_2$ sheet is zero outside of $0 \leq p \leq \dim X$, this will follow if we can show that 
\[ {a_{cdh}KH_q(-) \otimes \zll = 0} \quad  \textrm{ for all } \quad q < 0. \]
Recalling the isomorphism $K_n \otimes \zpi \cong KH_n \otimes \zpi$ mentioned at the beginning of the proof, it is equivalent to show that
\[ {a_{cdh}K_q(-) \otimes \zll = 0} \quad  \textrm{ for all } \quad q < 0. \]
This cdh sheaf is \ldh~separated (Lemma~\ref{lemm:KHfpslsep}) so the canonical morphism
\[ a_{cdh}K_q(-) \otimes \zll \to a_{ldh}K_q(-) \otimes \zll \]
is injective and it suffices to show that 
\[ a_{ldh}K_q(-) \otimes \zll = 0 \quad \textrm{ for each } \quad q < 0. \]
For any scheme $U$ in $Sch(X)$ and any section $s \in a_{ldh}K_q(U) \otimes \zll$ there exists an \ldh~covering $\{U_i \to U \}_{i \in I}$ such that each restriction $s|_{U_i}$ of the section $s$ is in the image of the canonical morphism of presheaves ${K_q(U_i) \otimes \zll \to a_{ldh}K_q (U_i)\otimes \zll}$. Each $U_i$ admits an \ldh~covering ${\{U_{ij} \to U_i\}_{j \in J_i}}$ with $U_{ij}$ regular (Corollary~\ref{coro:regularlCover}). Since $K_q$ vanishes on every regular scheme for $q < 0$ \cite[Proposition 6.8]{TT90}, the restrictions $s|_{U_{ij}}$ are all zero if $q < 0$, and therefore $s$ is also zero if $q < 0$.
\end{proof}


\bibliographystyle{amsalpha}
\bibliography{biblio}

\end{document}